\documentclass[]{article}

\usepackage{verbatim}
\usepackage{amssymb}
\usepackage{amsbsy}
\usepackage{amsmath}
\usepackage{amsthm}
\usepackage[mathscr]{eucal}
\expandafter\ifx\csname ifdraftMine\endcsname\relax
	\usepackage[dvipdfmx]{graphicx}
	\usepackage{caption}
\usepackage{subcaption}
	\usepackage{amsfonts}
  	\usepackage{enumitem}
	\usepackage{bm}
	\usepackage{listliketab}
\usepackage{bm}

\numberwithin{defn}{section}
\newcommand{\x}{\bm{x}}
\newcommand{\y}{\bm{y}}

\newcommand{\n}{\bm{n}}

\title{On the inclusion of damping terms in\\ the hyperbolic MBO algorithm}

\author{Elliott Ginder$^{*}$ (Meiji University)\\ Ayumu Katayama (Hokkaido University)}
\date{}
\begin{document}

%\begin{abstract}
%The hyperbolic MBO is a threshold dynamic algorithm which approximates interfacial motion by hyperbolic mean curvature flow. We introduce a generalization of this algorithm for imparting damping terms onto the equation of motion. We also construct corresponding numerical methods, and perform numerical tests. We also use our results to show that the generalized hyperbolic MBO is able to approximate motion by the standard mean curvature flow.
%\end{abstract}

\maketitle

\section{Introduction}
\fi
We develop an approximation method for computing the damped motion of interfaces under hyperbolic mean curvature flow (HCMF):
\begin{align}\label{hmcf}
	\alpha \x_{tt}(t,s)+\beta\x_{t}(t,s) = -\gamma\kappa(t,s)\nu (t,s).
\end{align}
In the above, $\x: [0,T)\times I \rightarrow {\bf{R}}^2$ denotes a closed curve in ${\bf{R}}^2$ (parameterized over an interval $I$), $T>0$ is a final time, $\kappa$ denotes the curvature of the interface, and $\nu$ is the outward unit normal of the interface. The nonnegative parameters $\alpha,\beta,$ and $\gamma$, designate mass, damping, and surface tension coefficients, respectively. The subscripts signify differentiation with respect to their variables, so that $\x_{tt}$ refers to the normal acceleration of the interface, and $\x_t$ denotes the normal velocity. We remark that the presence of the inertial term signifies that the HMCF is an oscillatory interfacial motion. 

The equation of motion (\ref{hmcf}) is accompanied by two initial conditions: one for the initial shape of the interface, and another prescribing the initial velocity field along the interface. It can be shown \cite{lefloch} that, when the initial velocity field is normal to the interface, the velocity field of the interface remains normal for the remainder of the flow. Although tangental velocities can be used to impart features such as rotation into the interfacial dynamics, our study assumes the initial velocity field to act in the normal direction of the interface.

\section{A generalized HMBO algorithm}
The original threshold dynamical (TD) algorithm (the so-called MBO algorithm, see \cite{MBO}) is a method for approximating motion by mean curvature flow (MCF). Borrowing on such ideas, a TD algorithm for hyperbolic mean curvature flow was introduced in \cite{label:jjiam}. Whereas previous TD algorithms utilize properties of the diffusion equation to approximate MCF, properties of wave propagation (along with a particular choice of initial condition) were used to design an approximation method for HMCF. For a time step size $\tau>0$, the error of the approximation was shown to be of the order $O(\tau).$ In this study, we will use properties of wave propagation, together with a suitable initial velocity field, to incorporate damping terms into the HMCF.  

Let time be discretized with a step size $\tau>0$, and $n$ be a non-negative integer. For the sake of simplicity in the exposition, let $V_n$ denote the normal velocity of the interface at the time step $n$, $\dot{V}_n$ be the normal acceleration, and $\kappa_n$ be the corresponding curvature of the interface. For the time being, we take the mass, damping, and surface tension coefficients to be unity, and proceed to construct an approximation method for the following interfacial dynamics:
\begin{align}
	\dot{V}_{n}-V_{n}=-\kappa _{n} .
	\label{label:damp}
\end{align}

%第3節においてはLevel-set関数を用いてHMCFを導出したが，ここでは波動方程式の解の公式を用いて式\ (\ref{label:damp})を導出する．また，簡単のため2次元領域において計算するが，3次元領域においても同様に計算できる．まず，波動方程式の初期値問題を次のように設定する．
Our approach is to observe the propagation of interfaces under the wave equation:
\begin{align}
	\begin{cases}
		u_{tt}=c^{2}\Delta u,\ &\text{in}\ (0,\tau)\times \Omega \\
		u(0,\bm{x})=u_{0}(\bm{x}),\ &\text{in}\  \Omega \\
		u_{t}(0,\bm{x})=-v_{0}(\bm{x}),\ &\text{in}\  \Omega\\
		\partial_{\n}u =0  &\text{on}\ (0,\tau)\times\partial\Omega,
	\end{cases}
\end{align}
where $\Omega$ is a given domain with smooth boundary, $c^2$ sets the wave speed, $u_0$ is an initial profile, $v_0$ designates the initial velocity, and $\tau$ is the time step. Although we have prescribed a Neumann boundary condition, $\partial_{\n}u =0$, we will only focus on the motion of interfaces located away from the boundary of the domain. In particular, away from the boundary, the short-time solution of the wave equation can be expressed using the Poisson formula:
\begin{align}
	u(t,\bm{x})=\dfrac{1}{2\pi ct}\int _{B(\bm{x},ct)} \dfrac{u_{0}(\bm{y})+\nabla u_{0}(\bm{y})\cdot (\bm{y}-\bm{x})-tv_{0}(\bm{y})}{\sqrt{c^{2}t^{2}-|\bm{y}-\bm{x}|^{2}}}d\bm{y},
	\label{label:poisson}
\end{align}
where $B(\bm{x},ct)$ denotes the ball centered at $\x$ with radius $ct$. 

Let $\Gamma^{n}$ be the closed curve at time step $n$, described as the boundary of a set $S_n$, and denote its signed distance function by
\begin{align}
	d_n(\x)=\begin{cases}
		\inf_{\y\in \Gamma^n}||\x-\y|| \hspace{30pt} \x \in S_n\\
		-\inf_{\y\in \Gamma^n}||\x-\y|| \hspace{21pt}\text{otherwise.}
	\end{cases}
\end{align}

We remark that $d_0$ is constructed from the given initial configuration of the interface, and that $d_{-1}$ can be constructed using the initial velocity field along the interface. This allows us to define $u_0(\x)$ as follows, for any non-negative integer:
$$u_{0}(\bm{x})=2d_{n}(\bm{x})-d_{n-1}(\bm{x}).$$

\begin{figure}[htbp]
	\centering
	\includegraphics[scale=0.4]{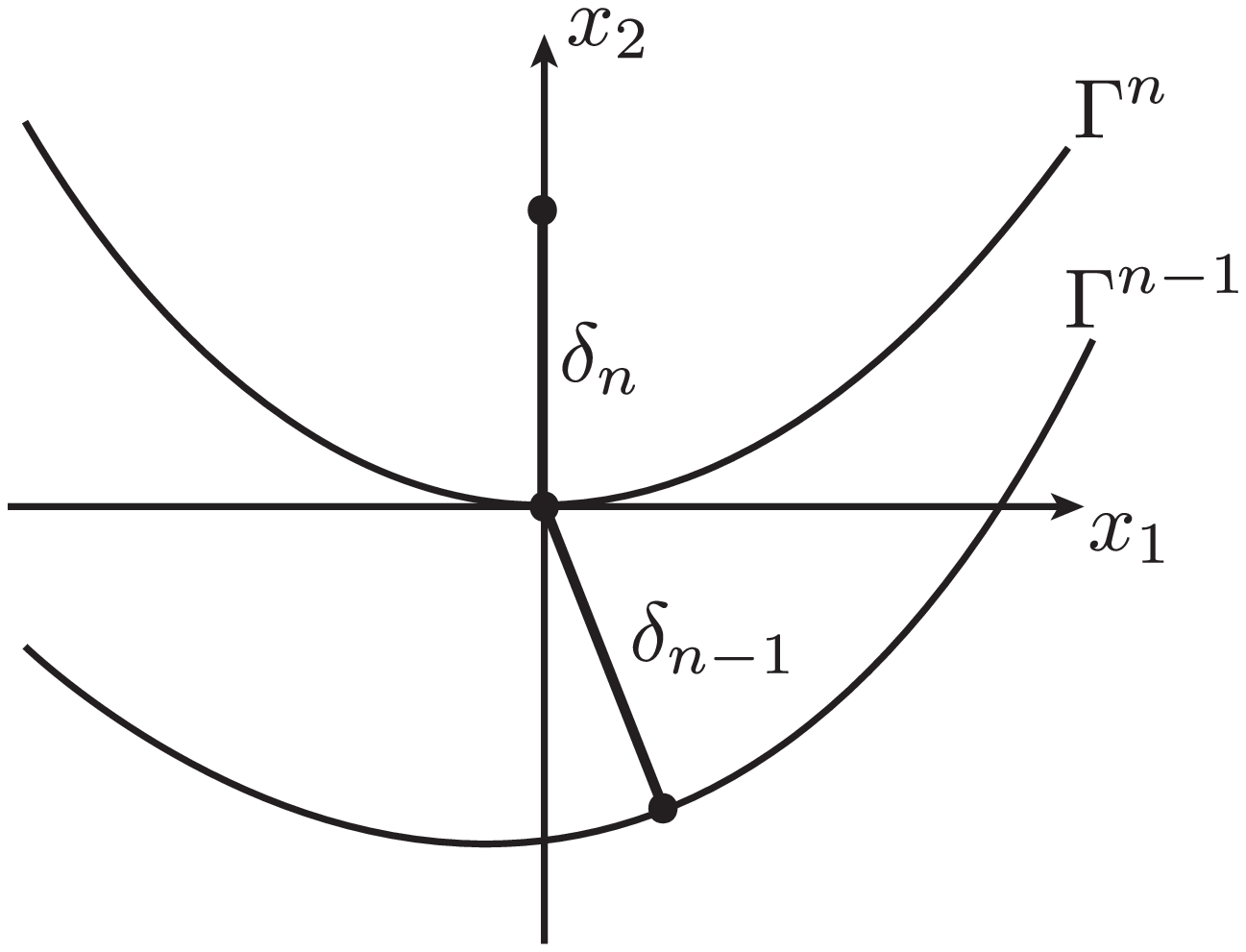}
	\caption{Motion of a single point of the interface in the normal direction．The point moves a distance $\delta_{n}$ at step $n.$ Without loss of generality, the direction of motion at the $n^{th}$ step is in the $x_2$ direction.}\label{fig:test}
\end{figure}

By taking $v_{0}(\bm{x})=0$ in (\ref{label:poisson}) and $c^2=2$, it can be shown (see \cite{label:jjiam}) that
\begin{align}
	\delta _{n}=\delta _{n-1}-(2\kappa _{n}-\kappa _{n-1})\tau ^{2}+O(\tau ^{3}),
	\label{label:deltan}
\end{align}
where $\delta_{n}$ denotes the distance traveled in the normal direction at time step $n$ (see figure \ref{fig:test}).
Denoting the average velocity within the time interval $[(n-1)\tau, n\tau]$ by $\bar{V}_{n}$, one has
\begin{align}
	\delta _{n}=\bar{V}_{n}\tau,\ \ \ \ \delta _{n-1}=\bar{V}_{n-1}\tau, \notag
\end{align}
and hence equation (\ref{label:deltan}) can be written:
\begin{align}
	\bar{V}_{n}\tau =\bar{V}_{n-1}\tau -(2\kappa _{n}-\kappa _{n-1})\tau ^{2}+O(\tau ^{3}).  \notag
\end{align}
Formally assuming $|\kappa _{n}-\kappa _{n-1}|<C\tau$ for some non-negative $C$, one obtains
\begin{align}
	\bar{V}_{n}\tau =\bar{V}_{n-1}\tau -\kappa _{n}\tau ^{2}+O(\tau ^{3}),
	\label{label:avev}
\end{align}
and hence
\begin{align}
	\dot{V}_{n} =-\kappa _{n}+O(\tau )\hspace{30pt}\text{(as $\tau \rightarrow 0$)}.  
	\label{label:approximate}
\end{align}

The damping term in equation (\ref{label:damp}) can be included by prescribing the initial velocity of the wave equation to be  $v_{0}(\bm{x})=d_{n}(\bm{x})$. This can be seen by expanding $d_{n}(\bm{x})$ in a Taylor series about $\bm{x}={\bf{0}}$ (see \cite{label:sdf}):
\begin{align}
	d_{n}(x_{1},x_{2})=x_{2}+\dfrac{1}{2}\kappa _{n}x_{1}^{2}+\dfrac{1}{6}(\kappa _{n})_{x_{1}}x_{1}^{3}-\dfrac{1}{2}\kappa _{n}^{2}x_{1}^{2}x_{2}+O(|\bm{x}|^{4}),
\end{align}
and appealing to Poisson's formula (\ref{label:poisson}):
\begin{align}
	u_{v}(t,\bm{x})=\dfrac{1}{2\pi ct}\int _{B(\bm{x},ct)} \dfrac{-tv_{0}(\bm{y})}{\sqrt{c^{2}t^{2}-|\bm{y}-\bm{x}|^{2}}}d\bm{y}.
\end{align}
Making the change of variables:
\begin{align}
	\bm{y}-\bm{x}&=ct\bm{z},
%	\Leftrightarrow 
%	\begin{pmatrix}
%		y_{1} \\
%		y_{2}
%	\end{pmatrix}
%	=&\begin{pmatrix}
%		ctz_{1}+x_{1} \\
%		ctz_{2}+x_{2}
%	\end{pmatrix}.
	\label{label:hh}
\end{align}
we note that
\begin{align}
	O(|\bm{y}|^{4})=O(t^{4})\hspace{30pt}\text{(as $t\rightarrow 0$)}.
	\label{label:order}
\end{align}
We thus investigate the contribution of the first four terms in the Taylor expansion, $u_{v}^{1},u_{v}^{2}, u_{v}^{3},$ and $u_{v}^{4}$. We begin with the lowest order term:
\begin{align}
	u_{v}^{1}(t,\bm{x})&=\dfrac{1}{2\pi ct}\int _{B(\bm{x},ct)} \dfrac{-ty_{2}}{\sqrt{c^{2}t^{2}-|\bm{y}-\bm{x}|^{2}}}d\bm{y} \notag \\
	&=\dfrac{-t}{2\pi ct}\int _{B(0,1)} \dfrac{ctz_{2}+x_{2}}{ct\sqrt{1-|\bm{z}|^{2}}}(ct)^{2}d\bm{z}. \notag
\end{align}
Appealing to function parity we have:
\begin{align*}
	\int _{B(0,1)} \dfrac{z_{2}}{\sqrt{1-|\bm{z}|^{2}}}d\bm{z}=0,
\end{align*}
and it follows that
\begin{align*}
	u_{v}^{1}(t,\bm{x})&=\dfrac{-t}{2\pi }\int _{B(0,1)} \dfrac{x_{2}}{\sqrt{1-|\bm{z}|^{2}}}d\bm{z}. 
\end{align*}
By making the change of variables:
\begin{align}
		z_{1}=r\text{cos}\theta, \hspace{15pt}
		z_{2}=r\text{sin}\theta,
	\label{label:varch1}
\end{align}
one arrives at
\begin{align}
	u_{v}^{1}(t,\bm{x})&=\dfrac{-tx_{2}}{2\pi }\int _{0}^{1}\int _{0}^{2\pi } \dfrac{1}{\sqrt{1-r^{2}}}rd\theta dr \notag \\
	&=\dfrac{-tx_{2}}{2\pi }\int _{0}^{1} \dfrac{2\pi r}{\sqrt{1-r^{2}}}dr \notag \\
	&=-tx_{2}\int _{0}^{1} \dfrac{r}{\sqrt{1-r^{2}}}dr. \notag 
\end{align}
Composite function integration yields
\begin{align}
	u_{v}^{1}(t,\bm{x})&=-tx_{2}\int _{0}^{1} (1-r^{2})^{-\frac{1}{2}}(-2r)\left( -\dfrac{1}{2} \right) dr
%	&=-tx_{2}\left( -\dfrac{1}{2} \right) 2(0-1). \notag \\
	=-tx_{2}.
	\label{label:v1}
\end{align}

We next consider the influence of the second term:
\begin{align}
	u_{v}^{2}(t,\bm{x})&=\dfrac{1}{2\pi ct}\int _{B(\bm{x},ct)} \dfrac{-\dfrac{1}{2}t\kappa _{n}y_{1}^{2}}{\sqrt{c^{2}t^{2}-|\bm{y}-\bm{x}|^{2}}}d\bm{y} \notag \\
	&=\dfrac{-t\kappa _{n}}{4\pi ct}\int _{B(\bm{x},ct)} \dfrac{y_{1}^{2}}{\sqrt{c^{2}t^{2}-|\bm{y}-\bm{x}|^{2}}}d\bm{y} \notag \\
	&=\dfrac{-t\kappa _{n}}{4\pi }\int _{B(0,1)} \dfrac{c^{2}t^{2}z_{1}^{2}+2ctx_{1}z_{1}+x_{1}^{2}}{\sqrt{1-|\bm{z}|^{2}}}d\bm{z}. \notag 
\end{align}
As before, function parity yields
\begin{align*}
	\int _{B(0,1)} \dfrac{z_{1}}{\sqrt{1-|\bm{z}|^{2}}}d\bm{z}=0,
\end{align*}
and hence
\begin{align*}
	u_{v}^{2}(t,\bm{x})&=\dfrac{-t\kappa _{n}}{4\pi }\int _{B(0,1)} \dfrac{c^{2}t^{2}z_{1}^{2}+x_{1}^{2}}{\sqrt{1-|\bm{z}|^{2}}}d\bm{z}. \notag 
\end{align*}
Making the change of variables (\ref{label:varch1}), we have
\begin{align}
	u_{v}^{2}(t,\bm{x})&=\dfrac{-t\kappa _{n}}{4\pi }\int _{0}^{1}\int _{0}^{2\pi } \dfrac{c^{2}t^{2}r^{3}\text{cos}^{2}\theta +x_{1}^{2}r}{\sqrt{1-r^{2}}}d\theta dr \notag \\
	&=\dfrac{-t\kappa _{n}}{4\pi }\int _{0}^{1}\int _{0}^{2\pi } \dfrac{c^{2}t^{2}r^{3}\dfrac{1+\text{cos}2\theta }{2} +x_{1}^{2}r}{\sqrt{1-r^{2}}}d\theta dr \notag \\
	&=\dfrac{-t\kappa _{n}}{4\pi }\int _{0}^{1} \dfrac{\pi \left(c^{2}t^{2}r^{3} +2x_{1}^{2}r\right)}{\sqrt{1-r^{2}}}dr. \notag 
\end{align}
Using another change of variables:
\begin{align}
	r=\text{cos}\theta,
	\label{label:varch2}
\end{align}
allows one to obtain:
\begin{align}
	u_{v}^{2}(t,\bm{x})&=\dfrac{-t\kappa _{n}}{4}\int _{0}^{\frac{\pi}{2}} \dfrac{c^{2}t^{2}\text{cos}^{3}\theta +2x_{1}^{2}\text{cos}\theta}{\sqrt{1-\text{cos}^{2}\theta}}\text{sin}\theta d\theta  \notag \\
	&=\dfrac{-t\kappa _{n}}{4}\int _{0}^{\frac{\pi}{2}} \left( c^{2}t^{2}\dfrac{\text{cos}3\theta +3\text{cos}\theta }{4} +2x_{1}^{2}\text{cos}\theta \right) d\theta  \notag \\
	&=\dfrac{-t\kappa _{n}}{4}\left( c^{2}t^{2}\left( -\dfrac{1}{12} +\dfrac{3}{4} \right) +2x_{1}^{2} \right) \notag \\
	&=-t\kappa _{n}\left( \dfrac{c^{2}t^{2}}{6} +\dfrac{x_{1}^{2}}{2}\right) .
	\label{label:v2}
\end{align}

The third term is similar:
\begin{align}
	u_{v}^{3}(t,\bm{x})&=\dfrac{1}{2\pi ct}\int _{B(\bm{x},ct)} \dfrac{-\dfrac{1}{6}t(\kappa _{n})_{x_{1}}y_{1}^{3}}{\sqrt{c^{2}t^{2}-|\bm{y}-\bm{x}|^{2}}}d\bm{y} \notag \\
	&=\dfrac{-t(\kappa _{n})_{x_{1}}}{12\pi ct}\int _{B(\bm{x},ct)} \dfrac{y_{1}^{3}}{\sqrt{c^{2}t^{2}-|\bm{y}-\bm{x}|^{2}}}d\bm{y} \notag \\
	&=\dfrac{-t(\kappa _{n})_{x_{1}}}{12\pi }\int _{B(0,1)} \dfrac{z_{1}^{3}+3c^{2}t^{2}x_{1}z_{1}^{2}+3ctx_{1}^{2}z_{1}+x_{1}^{3}}{\sqrt{1-|\bm{z}|^{2}}}d\bm{z}. \notag 
\end{align}
Again appealing to function parity:
\begin{align*}
	\int _{B(0,1)} \dfrac{z_{1}^{3}}{\sqrt{1-|\bm{z}|^{2}}}d\bm{z}=0,
\end{align*}
and therefore
\begin{align*}
	u_{v}^{3}(t,\bm{x})&=\dfrac{-t(\kappa _{n})_{x_{1}}}{12\pi }\int _{B(0,1)} \dfrac{3c^{2}t^{2}x_{1}z_{1}^{2}+x_{1}^{3}}{\sqrt{1-|\bm{z}|^{2}}}d\bm{z}. \notag
\end{align*}
The change of variables \ (\ref{label:varch1}) gives
\begin{align}
	u_{v}^{3}(t,\bm{x})&=\dfrac{-t(\kappa _{n})_{x_{1}}}{12\pi }\int _{0}^{1}\int _{0}^{2\pi } \dfrac{3c^{2}t^{2}x_{1}r^{3}\text{cos}^{2}\theta +x_{1}^{3}r}{\sqrt{1-r^{2}}}d\theta dr \notag \\
	&=\dfrac{-t(\kappa _{n})_{x_{1}}}{12\pi }\int _{0}^{1}\int _{0}^{2\pi } \dfrac{3c^{2}t^{2}x_{1}r^{3}\dfrac{1+\text{cos}2\theta }{2} +x_{1}^{3}r}{\sqrt{1-r^{2}}}d\theta dr \notag \\
	&=\dfrac{-t(\kappa _{n})_{x_{1}}}{12\pi }\int _{0}^{1} \dfrac{\pi \left(3c^{2}t^{2}x_{1}r^{3} +2x_{1}^{3}r\right)}{\sqrt{1-r^{2}}}dr, \notag
\end{align}
while (\ref{label:varch2}) allows one to express:
\begin{align}
	u_{v}^{3}(t,\bm{x})&=\dfrac{-t(\kappa _{n})_{x_{1}}}{12}\int _{0}^{\frac{\pi}{2}} \dfrac{3c^{2}t^{2}x_{1}\text{cos}^{3}\theta +2x_{1}^{3}\text{cos}\theta}{\sqrt{1-\text{cos}^{2}\theta}}\text{sin}\theta d\theta  \notag \\
	&=\dfrac{-t(\kappa _{n})_{x_{1}}}{12}\int _{0}^{\frac{\pi}{2}} \left( 3c^{2}t^{2}x_{1}\dfrac{\text{cos}3\theta +3\text{cos}\theta }{4} +2x_{1}^{3}\text{cos}\theta \right) d\theta  \notag \\
	&=\dfrac{-t(\kappa _{n})_{x_{1}}}{12}\left( 3c^{2}t^{2}x_{1}\left( -\dfrac{1}{12} +\dfrac{3}{4} \right) +2x_{1}^{3} \right) \notag \\
	&=\dfrac{-t(\kappa _{n})_{x_{1}}}{6}\left( c^{2}t^{2}x_{1} +x_{1}^{3}\right).
	\label{label:v3}
\end{align}

The final term follows the same approach:
\begin{align}
	u_{v}^{4}(t,\bm{x})&=\dfrac{1}{2\pi ct}\int _{B(\bm{x},ct)} \dfrac{\dfrac{1}{2}t\kappa _{n}^{2}y_{1}^{2}y_{2}}{\sqrt{c^{2}t^{2}-|\bm{y}-\bm{x}|^{2}}}d\bm{y} \notag \\
	&=\dfrac{t\kappa _{n}^{2}}{4\pi ct}\int _{B(\bm{x},ct)} \dfrac{y_{1}^{2}y_{2}}{\sqrt{c^{2}t^{2}-|\bm{y}-\bm{x}|^{2}}}d\bm{y} \notag \\
	&=\dfrac{t\kappa _{n}^{2}}{4\pi }\int _{B(0,1)} \dfrac{c^{3}t^{3}z_{1}^{2}z_{2}+c^{2}t^{2}x_{2}z_{1}^{2}+2c^{2}t^{2}x_{1}z_{1}z_{2}}{\sqrt{1-|\bm{z}|^{2}}}\notag\\&+\dfrac{2ctx_{1}x_{2}z_{1}+ctx_{1}^{2}z_{2}+x_{1}^{2}x_{2}}{\sqrt{1-|\bm{z}|^{2}}}d\bm{z} . \notag 
\end{align}
Function parity tells us that
\begin{align*}
	\int _{B(0,1)} \dfrac{z_{1}^{2}z_{2}}{\sqrt{1-|\bm{z}|^{2}}}d\bm{z}=0,
\end{align*}
and hence
\begin{align*}
	u_{v}^{4}(t,\bm{x})&=\dfrac{t\kappa _{n}^{2}}{4\pi }\int _{B(0,1)} \dfrac{c^{2}t^{2}x_{2}z_{1}^{2}+2c^{2}t^{2}x_{1}z_{1}z_{2}+x_{1}^{2}x_{2}}{\sqrt{1-|\bm{z}|^{2}}}d\bm{z}. \notag
\end{align*}
Applying the change of variables (\ref{label:varch1}) and computing gives
\begin{align}
	u_{v}^{4}(t,\bm{x})&=\dfrac{t\kappa _{n}^{2}}{4\pi }\int _{0}^{1}\int _{0}^{2\pi } \dfrac{c^{2}t^{2}x_{2}r^{3}\text{cos}^{2}\theta +2c^{2}t^{2}x_{1}r^{3}\text{cos}\theta \text{sin}\theta +x_{1}^{2}x_{2}r}{\sqrt{1-r^{2}}}d\theta dr \notag \\
	&=\dfrac{t\kappa _{n}^{2}}{4\pi }\int _{0}^{1}\int _{0}^{2\pi } \dfrac{c^{2}t^{2}x_{2}r^{3}\dfrac{1+\text{cos}2\theta }{2}+c^{2}t^{2}x_{1}r^{3}\text{sin}2\theta +x_{1}^{2}x_{2}r}{\sqrt{1-r^{2}}}d\theta dr \notag \\
	&=\dfrac{t\kappa _{n}^{2}}{4\pi }\int _{0}^{1} \dfrac{\pi \left(c^{2}t^{2}x_{2}r^{3} +2x_{1}^{2}x_{2}r\right)}{\sqrt{1-r^{2}}}dr. \notag
\end{align}
Using the change of variables (\ref{label:varch2}) leads us to the expession:
\begin{align}
	&=\dfrac{t\kappa _{n}^{2}}{4}\int _{0}^{\frac{\pi}{2}} \dfrac{c^{2}t^{2}x_{2}\text{cos}^{3}\theta +2x_{1}^{2}x_{2}\text{cos}\theta}{\sqrt{1-\text{cos}^{2}\theta}}\text{sin}\theta d\theta  \notag \\
	&=\dfrac{t\kappa _{n}^{2}}{4}\int _{0}^{\frac{\pi}{2}} \left( c^{2}t^{2}x_{2}\dfrac{\text{cos}3\theta +3\text{cos}\theta }{4} +2x_{1}^{2}x_{2}\text{cos}\theta \right) d\theta  \notag \\
	&=\dfrac{t\kappa _{n}^{2}}{4}\left( c^{2}t^{2}x_{2}\left( -\dfrac{1}{12} +\dfrac{3}{4} \right) +2x_{1}^{2}x_{2} \right) \notag \\
	&=t\kappa _{n}^{2}\left( \dfrac{c^{2}t^{2}x_{2}}{6} +\dfrac{x_{1}^{2}x_{2}}{2}\right).
	\label{label:v4}
\end{align}

Equations  (\ref{label:v1})，(\ref{label:v2})，(\ref{label:v3})，and (\ref{label:v4}) together express
\begin{align}
	u_{v}(t,\bm{x})=&-t\left( x_{2}+\kappa _{n}\left( \dfrac{c^{2}t^{2}}{6} +\dfrac{x_{1}^{2}}{2}\right)+\dfrac{(\kappa _{n})_{x_{1}}}{6}\left( c^{2}t^{2}x_{1} +x_{1}^{3}\right) \right)\\ &+t \kappa _{n}^{2}\left( \dfrac{c^{2}t^{2}x_{2}}{6} +\dfrac{x_{1}^{2}x_{2}}{2}\right). \notag
\end{align}
Upon taking $t=\tau$ and $\bm{x}=(0,\delta _{n})$, we arrive at
\begin{align}
	0&=-\tau \left( \delta _{n}+\kappa _{n}\dfrac{c^{2}\tau ^{2}}{6}-\kappa _{n}^{2}\dfrac{c^{2}\tau ^{2}\delta _{n}}{6} \right)
	  =-\delta _{n}\tau +O(\tau ^{3}). \label{label:uv} 
\end{align}
Combining this equation with our previous results yields: 
\begin{align}
	%&0=\delta _{n}-\delta _{n-1}+(2\kappa _{n}-\kappa _{n-1})\tau ^{2}-\delta _{n}\tau +O(\tau ^{3}), \notag \\
	%\Leftrightarrow \
	&\delta _{n}-\delta _{n-1}-\delta _{n}\tau =-(2\kappa _{n}-\kappa _{n-1})\tau ^{2}+O(\tau ^{3}).
	\label{label:delv}
\end{align}
Writing $\delta _{n}=\bar{V}_{n}\tau $ and $\delta _{n-1}=\bar{V}_{n-1}\tau$ in equation (\ref{label:delv}) expresses
\begin{align}
	\bar{V}_{n}\tau -\bar{V}_{n-1}\tau -\bar{V}_{n}\tau ^{2}=-(2\kappa _{n}-\kappa _{n-1})\tau ^{2}+O(\tau ^{3}). \notag
\end{align}
Formally assuming $|\kappa _{n}-\kappa _{n-1}|<C\tau$, for some constant $C$, and dividing both sides by $\tau ^{2}$ gives
\begin{align}
	\dfrac{\bar{V}_{n}-\bar{V}_{n-1}}{\tau } -\bar{V}_{n}=-\kappa _{n}+O(\tau ).  \notag
\end{align}
It follows that the damping term enters the equation of motion:
\begin{align}
	\dot{V}_{n} -V_{n}=-\kappa _{n}+O(\tau ). 
	\label{label:dampe}
\end{align}
By linearity, taking $u_{0}(\bm{x})=a(2d_{n}(\bm{x})-d_{n-1}(\bm{x}))$ and $v_{0}(\bm{x})=bd_{n}(\bm{x})$ one can obtain the interfacial
motion:
\begin{align}
	a\dot{V}_{n} -bV_{n}=-\dfrac{ac^{2}}{2}\kappa _{n}+O(\tau ),
	\label{label:dampp}
\end{align}
where $a$ and $b$ are real parameters. Therefore, one can rewrite the parameters:
\begin{align}
		\alpha =a, \hspace{15pt}
		\beta =-b, \hspace{15pt}
		\gamma =\dfrac{ac^{2}}{2},
	\label{label:para} 
\end{align}
to approximate a prescribed interfacial motion:
\begin{align}
	\alpha\dot{V} +\beta V=-\gamma \kappa+O(\tau ). 
	\label{label:gdamp}
\end{align}

The previous results show that the wave equation's initial velocity can be used in the HMBO algorithm to impart damping terms. In the next section, by choosing parameters, we will make a numerical investigation into using the HMBO to approximate interfacial motion by the standard mean curvature flow.

\section{The HMBO approximation of mean curvature flow}
%本項ではHMBOアルゴリズムを用いたMCFの数値計算について考える．運動方程式の導出については前項と同様に行い，数値計算結果も示す．
%\subsubsection{運動方程式の導出}
%MCFの運動方程式は式\ (\ref{label:mcf})で示している．これをHMBOアルゴリズムを用いて計算することを考える．HMCFの運動方程式について，一般的な形は式\ (\ref{label:gdamp})として得られているが，このとき$\tilde{a}=0$とすると，式\ (\ref{label:para})より同時に$\tilde{c}=0$となり，
%\begin{align}
%	\tilde{b}V_{n}=O(\tau ). \notag
%\end{align}
%となり，求める結果を得られない．
An approximation method for mean curvature flow can be obtained by returning to 
equation (\ref{label:poisson}) and choosing appropriate initial conditions. For a predetermined time
step $\tau>0$, we take $u_{0}(\bm{x})=0$, $v_{0}(\bm{x})=d_{n}(\bm{x})$, and $c^{2}= \lambda / \tau $. Then 
equation (\ref{label:uv}) gives
\begin{align}
	0&=-\tau \left( \delta _{n}+\kappa _{n}\dfrac{\lambda \tau}{6}-\kappa _{n}^{2}\dfrac{\lambda \delta _{n}\tau}{6} \right)+O(\tau ^{3}), \notag \\
	&=\delta _{n}+\kappa _{n}\dfrac{\lambda \tau}{6}-\kappa _{n}^{2}\dfrac{\lambda \delta _{n}\tau}{6}+O(\tau ^{2}). \notag
\end{align}
%このとき誤差は式\ (\ref{label:hh})，(\ref{label:order})から求められる．$\delta _{n}=\bar{V_{n}}\tau$とすると，
%\begin{align}
%	0&=\bar{V}_{n}\tau+\kappa _{n}\dfrac{\lambda \tau}{6}-\kappa _{n}^{2}\dfrac{\lambda \bar{V}_{n}\tau ^{2}}{6}+O(\tau ^{2}). \notag
%\end{align}
%従って，次式を得る．
Preceeding as in the previous section, we obtain
\begin{align}
	\bar{V}_{n}=-\dfrac{\lambda }{6}\kappa _{n}+O(\tau ). \notag
\end{align}
Since $\lambda $ is a free parameter, we find that the corresponding threshold dynamics can approximate
curvature flow with a parameter $\gamma$:
\begin{align}
	V_{n}=-\gamma\kappa _{n}+O(\tau )\hspace{30pt}(\text{as $\tau \rightarrow 0$}).  
	\label{label:mboh}
\end{align}
\section{Numerical investigation}
We will now perform a numerical error analysis of the HMBO approximation of MCF. The numerical method's performance will be compared to the case of a circle evolving by MCF. In such a setting, the evolution of the circle's radius is governed by the solution of the following ordinary differential equation:
%相の数を$N=2$とする．このとき2次元領域においてHMCFに従う運動について，その誤差評価は既に行われている\cite{label:jjiam}．ここでは2次元領域におけるHMBOアルゴリズムを用いたMCFの計算について誤差評価する．そのための例題を次式に示す．
\begin{align}
	\begin{cases}
		\dot{r}(t)=-\dfrac{1}{r(t)} \ \ \ \ t>0, \\
		r(0)=r_{0},
		\label{label:exr}
	\end{cases}
\end{align}
where $r_0$ is the initial radius of the circle. We remark that the radius decreases until its extinction time $t_{e}= r_{0}^{2}/2$, and that $r(t)= \sqrt{r_{0}^{2}-2t}$.
%計算領域が2次元なので，平均曲率$\kappa =1/r(t)$となっている．この問題は厳密解が知られているため誤差の評価が可能である．
%この微分方程式の厳密解を示す．これは単純な常微分方程式を解くだけでよい．
%\begin{align}
%	&\dot{r}(t)=-\dfrac{1}{r(t)}. \notag \\
%	\Rightarrow \ &\int r\dot{r}dr=\int -dt. \notag \\
%	\Rightarrow \ &\dfrac{r^{2}}{2}=-t+C. 
%	\label{label:cre}
%\end{align}
%$r(0)=r_{0}$なので，
%\begin{align*}
%	C=\dfrac{r_{0}^{2}}{2}.
%\end{align*}
%よって式\ (\ref{label:cre})より，厳密解は，
%\begin{align}
%	&\dfrac{r^{2}}{2}=\dfrac{r_{0}^{2}}{2}-t. \notag \\
%	\therefore \ &r= \sqrt{r_{0}^{2}-2t}. 
%	\label{label:re}
%\end{align}
%加えて，ここでは界面が消滅する時間についても考える．これは$r(t)=0$を満たすような$t$を計算することに他ならない．従って，界面が消滅する時間を$t_{e}$とすると，
%\begin{align}
%	t_{e}=\dfrac{r_{0}^{2}}{2}. 
%	\label{label:te}
%\end{align}
%本節では式\ (\ref{label:re})，(\ref{label:te})を用いて誤差を評価する．

The HMBO approximation method solves the following wave equation for a small time $\tau>0$:
\begin{align}
	\begin{cases}
		u_{tt}=c^{2}\Delta u\ &\text{in}\ (0,\tau)\times \Omega \\
		u(0,\bm{x})=0\ &\text{in}\  \Omega \\
		u_{t}(0,\bm{x})=d_{k}(\bm{x})\ & \text{in}\ \Omega \\
		\partial_{\boldsymbol{n}}u=0 \ &\text{on}\ \partial \Omega,
	\end{cases}
	\label{label:ehmbo}
\end{align}
where $\Omega =(-2,2)\times (-2,2)$ and $k$ denotes the $k^{th}$ step of the HMBO algorithm.
%界面の初期形状を設定する．相の領域を次のようにとる．
%\begin{align*}
%	&P_{1}^{0}=\{ \bm{x}\in \Omega \ |\ ||\bm{x}|| > 1 \}, \\
%	&P_{2}^{0}=\Omega \setminus P_{1}^{0}.
%\end{align*}
%つまり，界面の初期形状は原点を中心とする半径$r_{0}=1$の円周となり，符号付距離$d_{i}^{0}(\bm{x})$は，
We choose the initial interface to be a circle with radius one, so that 
\begin{align*}
	&d_{0}(\bm{x})=||\bm{x}||-1.
\end{align*}
The initial velocity at the $k^{th}$ step is then defined as the signed distance function to the zero level set of the solution to the wave equation: 
\begin{align}
	d_k(\x)=\begin{cases}
		\inf_{\y\in \partial\{u(\x,\tau)>0\}}||\x-\y|| \hspace{30pt} \x \in \{u(\x,\tau)>0\}\\
		-\inf_{\y\in \partial\{u(\x,\tau)>0\}}||\x-\y|| \hspace{21pt}\text{otherwise.}
	\end{cases}
\end{align}
%このとき補間パラメータは$\epsilon \geq 2\sqrt{2}$とし，領域全体で符号付距離を計算する．さらに，対象としている領域が2次元なので，界面は有限個の節点と線分から構成され，この線分との距離を計算することとなる．
%\subsection{空間分割数について}
%空間の分割数によって誤差がどの程度変化するか調べる．そのため計算の際に用いるパラメータを設定する．
Since the extinction time $t_e$ depends on $r_0$, we set $\tau=t_{e}/N_{\tau}$. 
%まず時間ステップ$\tau $について，界面が消滅する時刻$t_{e}$が界面の初期半径$r_{0}$に依存して決まるので，$\tau =t_{e}/N_{\tau}$としてある程度の精度を確保する．
Here $r_{0}=1$ (hence $t_{e}=0.5$), and we set $N_{\tau}=150$ to ensure a level of precision. The time step is then $\tau =3.33\times 10^{-3}$. The target problem (\ref{label:exr}) corresponds to $\gamma=1$ in equation (\ref{label:mboh}), and we thus set $c^{2}=6/\tau $. 

\begin{figure}[htbp]
	\centering
	\includegraphics[scale=0.25]{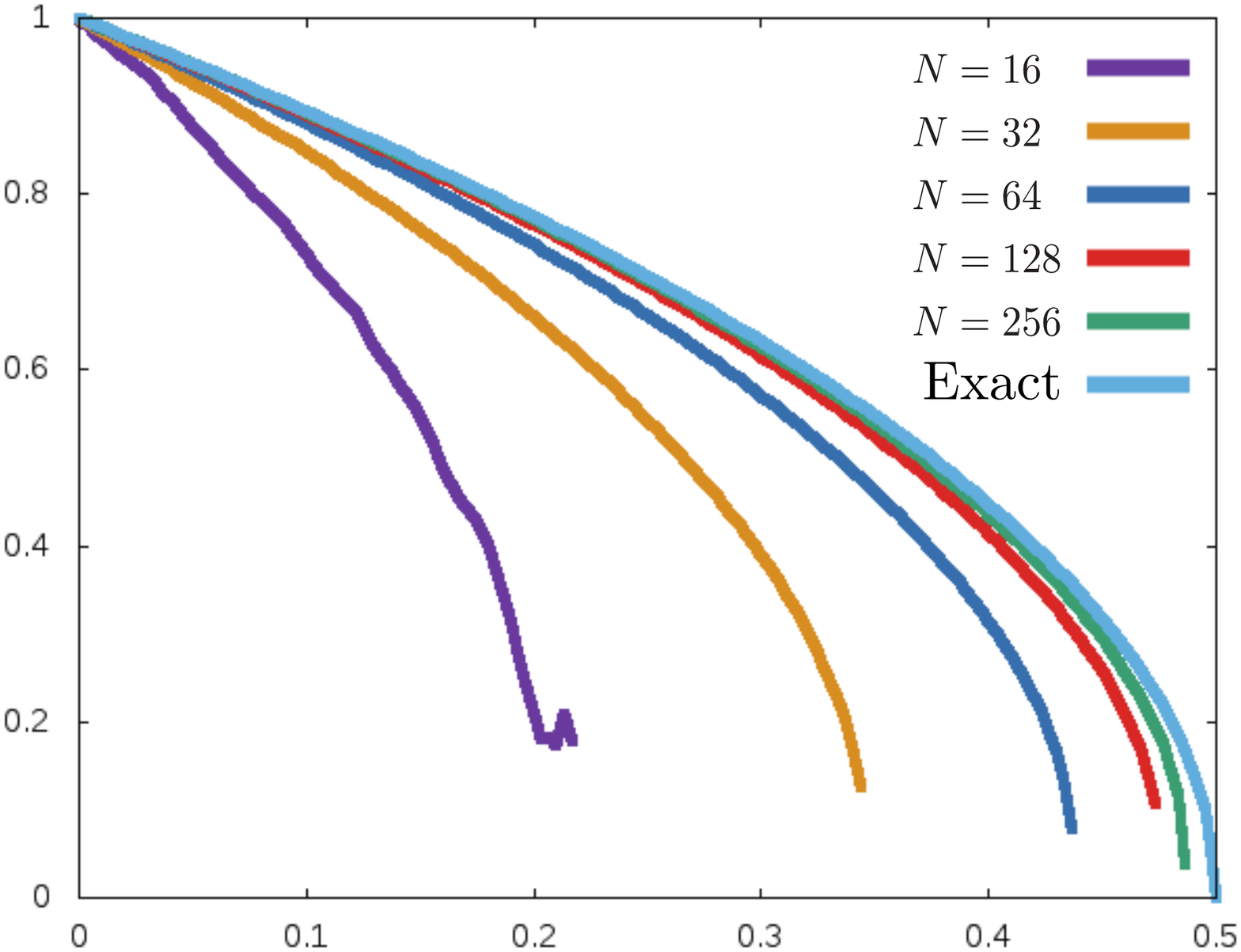}
	\caption{Convergence of the approximation method as $N$ is increased.}
	\label{label:figte}
\end{figure}

Finite differences are used to numerically solve the wave equation with a time step $\Delta t=2.22\times 10^{-6}$. The grid spacing in the $x$ and $y$ directions are equal to $\Delta x=2/(N-1)$, where $N$ is a natural number. We examine the numerical error when $N=2^{j},$ for $j=4,5,6,7,8$. The numerical results are shown in figure \ref{label:figte}, where the radius of the numerical solution is defined to be the average distance $\tilde{r}(t)$ of the level set's point cloud to the origin. 
%図\ (\ref{label:figte})から読み取れる2つの事実について考察する．まず1つ目に分割数が小さいときは非常に精度が悪いということである．特に$N=16$の場合では数値の振動が見られる．分割数が小さい場合は界面の近似に用いられる線分が少なく，そのため精度が悪化する為である．さらにこれは分割数を大きくとった場合にも同じ状況が考えられる．界面の半径が小さくなると相対的に領域の分割数が小さくなり，界面の近似精度が徐々に低下していく．これは界面が消滅する際だけでなく，界面の曲率が大きい場合にも同様の状況が考えられる．
%次に，計算不可能となる界面の半径を分割数から直接読み取ることができないことを挙げる．図\ (\ref{label:figte})に示した計算結果は界面が消滅する1ステップ前までのデータをグラフ化しているが，$N=64$と$N=128$のグラフを比較すると分割数が大きい場合の方がより大きな半径で計算が終了している．ただし厳密解との誤差は減少しているので，大きな問題とはならない．
The error is measured using the quantity:
\begin{align}
	Err(t)=\int _{0}^{T}| r(t)-\tilde{r}(t)|dt.
	\label{label:err}
\end{align}
Since the extinction time of the numerical solution differs from the exact solution, the actual error is computed as follows:
\begin{align}
	Err(t)\approx \displaystyle\sum_{i=0}^{N_{s}}|r(i\tau )-\tilde{r}(i\tau )|\tau,
	\label{label:erra}
\end{align}
where $N_{s}$ denotes the number of time steps until the numerical solution's radius disappears (the corresponding time is $N_{s}\tau$).
Our results are summarized in table (\ref{label:erra}), where we observe the convergence of our method to the exact solution.
\begin{table}[htb]
	\centering
	\caption{Error Table with respect to $\Delta x.$}
	\label{label:taberr1}
	\begin{tabular}{c|c|c}
		$N$ & $N_{s}\tau $ & $Err$ \\ \hline
%    		8 & 0.083333 & 0.014775 \\
		16 & 0.223333 & 0.044613\\
		32 & 0.343333 & 0.039463\\
		64 & 0.436667 & 0.022746\\
		128 & 0.473333 & 0.008509\\
		256 & 0.486667 & 0.003907\\
  	\end{tabular}
\end{table}

\section{Acknowledgments}
E. Ginder would like to acknowledge the support of JSPS Kakenhi Grant Number 17K14229, as well as that from the Presto Research Program of the Japan Science and Technology Agency.

\expandafter\ifx\csname ifdraftMine\endcsname\relax
  \end{document}
\fi

%%%%%%%%%%%%%%%%%%%%%%%%%%%%%%%%%
% References
%%%%%%%%%%%%%%%%%%%%%%%%%%%%%%%%%